\documentclass{article}
\usepackage[english]{babel}
\usepackage{latexsym,amssymb}
\usepackage[cp1251]{inputenc}
\usepackage{amsfonts,amssymb}
\usepackage{euscript}
\usepackage[multiple]{footmisc}
\newcounter{theorem} %[section]
\newcounter{lemma} %[section]

\begin{document}

\vspace*{7mm}

\Large

%\addtocounter{page}{-1}
%\renewcommand{\baselinestretch}{1.35}
\thispagestyle{empty}

%\noindent УДК 517.548\\
\Large
 \begin{center}
    \textbf{Vitalii Shpakivskyi\footnote{Institute of Mathematics of the National Academy of Sciences of Ukraine}\footnote{Zhytomyr Polytechnic State University}}
       \end{center}
  %  \vspace{2mm}
     \begin{center}
    shpakivskyi86@gmail.com
    \end{center}
%\vspace{3mm}
\vspace{7mm}
\Huge
 \begin{center}
   \textbf{
PART 4}
 \end{center}
\vspace{3mm}
%\vspace{7mm}

\LARGE
 \begin{center}
   \textbf{Hypercomplex Method for Solving Nonlinear PDEs}
 \end{center}
\vspace{5mm}

% У роботі розвинуто метод П.~Гумберта  розв'язування диференціальних рівнянь з частинними похідними. Даний метод застосовано до
% побудови розв'язків тривимірних рівнянь Лапласа, Гельмгольца і Пуассона в класах голоморфних функцій.

%\vspace{15mm}

%\noindent {\bf  Vitalii Shpakivskyi}
%\vspace{3mm}

%\noindent {\bf Holomorphic functions of several hypercomplex variables and their applications}

%\vspace{7mm}

% The work developed P.~Humbert's method of solving PDEs. This method is applied to
% construction of solutions of the three-dimensional equations of Laplace, of Helmholtz and of Poisson in the classes of holomorphic functions.
 
\Large
%\newpage 
%\section{Основні поняття}

\section{Algebras $\mathcal{A}_2\bigl(u(x,y)\bigr)$}

In Part 2, we considered algebras with variable structural constants, that is, $A$ and $B$ were functions of $x,y$, and proposed applications of these algebras to linear operator-differential PDEs with variable coefficients. In this part, we introduce a two-dimensional commutative associative algebra with structural constants $A$ and $B$, where the ``constants'' $A$ and $B$ depend not on the independent variables $x,y$, but on the value of a certain function $u(x,y)$. Monogenic functions will be introduced in this algebra, and applications of such monogenic functions to the construction of solutions to a certain class of nonlinear PDEs will be proposed.

Let $D$ be a domain in $\mathbb{R}^2$ and let $u:D\to I$ be a function, where $I\subseteq\mathbb{R}$ is a finite or infinite interval (domain). Suppose that $u\in C^2(D)$ and that $A,B\in C^{1}(I)$. Then, to each point $(x,y)$, we associate a two-dimensional commutative associative algebra
$\mathcal{A}_2\bigl(u(x,y)\bigr)$ with basis $\{e_1,e_2\}$ and multiplication table
\begin{equation}\label{tabl-mnozh-2-vym-algebra-neliniyni}
\begin{array}{c||c|c|}
\cdot & e_1 & e_2\\
\hline\hline
e_1 & e_1 & e_2\\
\hline
e_2 & e_2 & A(u)e_1+B(u)e_2\\
\hline
\end{array}\,.
\end{equation}

The algebra $\mathcal{A}_2\bigl(u(x,y)\bigr)$   is a generalization of the classical concept of an algebra \cite{Pierce}.

\section{Monogenic functions in the commutative algebra
$\mathcal{A}_2\bigl(u(x,y)\bigr)$}

Next, we generalize elements of the theory of monogenic functions developed in the monograph \cite{Plaksa-Shpakivskyi-2023} to the algebra $\mathcal{A}_2\bigl(u(x,y)\bigr)$.

In what follows, for simplicity of presentation and understanding, we consider the variable
$\zeta=xe_1+ye_2$, $(x,y)\in D$, where $\{e_1,e_2\}$ is the common fixed basis of all algebras in the family
$\mathcal{A}_2\bigl(u(x,y)\bigr)$.

Consider a function $\Phi:D\to \mathcal{A}_2\bigl(u\bigr)$.
Since all algebras in the family $\mathcal{A}_2\bigl(u\bigr)$ have the common basis $\{e_1,e_2\}$, the function $\Phi$ can be represented in the form
\begin{equation}\label{rozklad za basysom -2-vym-neliniy}
\Phi(\zeta)=u(x,y)e_1+v(x,y)e_2,
\end{equation}
where the first component $u(x,y)$ is the same function as that appearing in the multiplication table of the algebra
(\ref{tabl-mnozh-2-vym-algebra-neliniyni}); $u:D\to I$, $v:D\to\mathbb{R}$.
The functions $u$ and $v$ will be called the components of the function $\Phi$.

%\textbf{Important remark.} In the representation
%(\ref{rozklad za basysom -2-vym-neliniy}), the first component $u(x,y)$
%is exactly the same function as that appearing in the multiplication table
%(\ref{tabl-mnozh-2-vym-algebra-neliniyni}).

\vskip2mm
\textbf{Definition 1.}
Let the components $u,v$ of the function
(\ref{rozklad za basysom -2-vym-neliniy}) belong to $C^2(D)$.
The function $\Phi$ will be called monogenic in the domain $D$ if, at every point $(x,y)\in D$, there exists an element
$\Phi'(\zeta)\in \mathcal{A}_2\bigl(u\bigr)$ such that
\begin{equation}\label{ozn-monog-A-2-zmin-koef-nelin}
d\Phi=\Phi'(\zeta)d\zeta.
\end{equation}
The element $\Phi'(\zeta)$ will be called the derivative of the function $\Phi$ at the point $\zeta$.

We have
$$
d\Phi=\Phi_xdx+\Phi_ydy, \qquad
d\zeta=dx e_1+dy e_2.
$$
Therefore, from equality
(\ref{ozn-monog-A-2-zmin-koef-nelin}) we obtain
$$
\Phi_xdx+\Phi_ydy
=
\Phi'(\zeta)dx e_1+\Phi'(\zeta)dy e_2.
$$
Since $e_1$ is the unity of the algebra, comparison of the coefficients of $dx$ gives
$\Phi'(\zeta)=\Phi_x$.
Substituting this expression into the equality of the coefficients of $dy$, we obtain
$\Phi_y=\Phi_xe_2$.

Thus, we have the following theorem.

\vskip2mm
\textbf{Theorem 1.}
A function $\Phi$ of the form
(\ref{rozklad za basysom -2-vym-neliniy}) is monogenic in the domain $D$ if and only if
\begin{equation}\label{A-2-zmin-koef-umovy-K-R-neliniy}
\Phi_y=\Phi_xe_2.
\end{equation}
Conditions (\ref{A-2-zmin-koef-umovy-K-R-neliniy}) are analogues of the Cauchy--Riemann conditions.

Let us write the Cauchy--Riemann conditions in terms of the components $u,v$ of a monogenic function.

Indeed,
$$
u_ye_1+v_ye_2
=
A\bigl(u(x,y)\bigr)v_xe_1
+
\bigl(u_x+B\bigl(u(x,y)\bigr)v_x\bigr)e_2.
$$
Hence, equating the coefficients of $e_1$ and $e_2$, we obtain the system
\begin{equation}\label{A-2-zmin-koef-systema-K-R-nelin}
u_y=A\bigl(u(x,y)\bigr)v_x,
\qquad
v_y=u_x+B\bigl(u(x,y)\bigr)v_x.
\end{equation}

Clearly, the Cauchy--Riemann conditions
(\ref{A-2-zmin-koef-umovy-K-R-neliniy}) or
(\ref{A-2-zmin-koef-systema-K-R-nelin})
depend on the value of the first component of the monogenic function.
Therefore, even within the same domain $D$, their coefficients vary from point to point.

\section{Monogenic functions and nonlinear PDEs}

Let us establish which differential equations are satisfied by the components of the monogenic function (\ref{ozn-monog-A-2-zmin-koef-nelin}).

Differentiating the first equation of system
(\ref{A-2-zmin-koef-systema-K-R-nelin})
with respect to the variables $x$ and $y$, we obtain
\begin{equation}\label{A-2-zmin-koef-systema-K-R-nelin-rozsh-1}
u_{xy}=A'(u)u_xv_x+A(u)v_{xx}
\end{equation}
and
\begin{equation}\label{A-2-zmin-koef-systema-K-R-nelin-rozsh-2}
u_{yy}=A'(u)u_yv_x+A(u)v_{xy}.
\end{equation}

Differentiating the second equation of system
(\ref{A-2-zmin-koef-systema-K-R-nelin})
with respect to $x$, we have
\begin{equation}\label{A-2-zmin-koef-systema-K-R-nelin-rozsh-3}
v_{xy}
=
u_{xx}+B'(u)u_xv_x+B(u)v_{xx}.
\end{equation}

Substituting (\ref{A-2-zmin-koef-systema-K-R-nelin-rozsh-3}) into
(\ref{A-2-zmin-koef-systema-K-R-nelin-rozsh-2}), we obtain
$$
u_{yy}
=
A(u)u_{xx}
+
\bigl(A'(u)u_y+A(u)B'(u)u_x\bigr)v_x
+
A(u)B(u)v_{xx}.
$$
Multiplying this equality by $A(u)$, and using the first equation of system
(\ref{A-2-zmin-koef-systema-K-R-nelin}) together with
(\ref{A-2-zmin-koef-systema-K-R-nelin-rozsh-1}), we have
$A(u)v_x=u_y$ and
$$
A(u)v_{xx}
=
u_{xy}-A'(u)u_xv_x.
$$
Therefore,
$$
A(u)u_{yy}
=
A^2(u)u_{xx}
+
A'(u)u_y^2
+
A(u)B'(u)u_xu_y
+
A(u)B(u)u_{xy}
-
B(u)A'(u)u_xu_y.
$$
Thus, the first component of a monogenic function satisfies the nonlinear differential equation
\begin{equation}\label{rivn-dlia-u-A-B-vid-u-neliniyne}
A(u)u_{yy}
=
A^2(u)u_{xx}
+
A(u)B(u)u_{xy}
+
A'(u)u_y^2
+
\Biggr(
A(u)B'(u)-B(u)A'(u)
\Biggr)u_xu_y.
\end{equation}
Here $A$ and $B$ are arbitrary functions of class $C^1(I)$.

Thus, the following theorem holds.

\vskip2mm
\textbf{Theorem 2.}
Let $A,B\in C^1(I)$, let the components $u,v\in C^2(D)$, and let the function
(\ref{rozklad za basysom -2-vym-neliniy})
be monogenic in the domain $D$. Then the first component $u$ of the function $\Phi$ satisfies the nonlinear differential equation
(\ref{rivn-dlia-u-A-B-vid-u-neliniyne}) in $D$.

\vskip2mm
\textbf{Remark.}
If, in addition, $A(u)\neq 0$ in the domain $D$, then equation
(\ref{rivn-dlia-u-A-B-vid-u-neliniyne}) can be divided by $A(u)$ and written in the normalized form
\begin{equation}\label{rivn-dlia-u-A-B-vid-u-neliniyne-normal}
u_{yy}
=
A(u)u_{xx}
+
B(u)u_{xy}
+
\frac{A'(u)}{A(u)}\,u_y^2
+
\left(
B'(u)
-
\frac{B(u)A'(u)}{A(u)}
\right)
u_xu_y.
\end{equation}

Equation (\ref{rivn-dlia-u-A-B-vid-u-neliniyne}) is a nonlinear second-order differential equation with variable coefficients depending on the unknown solution itself.
It will be shown below that the Cauchy--Riemann system
(\ref{A-2-zmin-koef-systema-K-R-nelin})
admits a hodograph transformation that reduces it to a linear system.

\section{Hodograph transformation}

The hodograph transformation is a classical method for reducing certain nonlinear partial differential equations and systems to linear or simpler equations;
 see, for example, \cite{Ames,Clarkson,Qu}.

We show that the nonlinear Cauchy--Riemann system
(\ref{A-2-zmin-koef-systema-K-R-nelin})
can be linearized by means of a hodograph transformation.

Suppose that the Jacobian
$$
J=\frac{\partial(u,v)}{\partial(x,y)}=u_xv_y-u_yv_x
$$
does not vanish in the domain $D$. Then, by the inverse function theorem, there exists a local inverse mapping
\begin{equation}\label{hodograf-nelin}
x=x(u,v),\qquad
y=y(u,v).
\end{equation}

\vskip2mm
\textbf{Lemma 1.}\label{Lema-1-nelin}
The Jacobian matrix of the inverse mapping is the inverse of the Jacobian matrix of the original mapping.

\textbf{Proof.}
Let
$$
F(x,y)=\bigl(u(x,y),\,v(x,y)\bigr),
$$
and let
$$
G(u,v)=\bigl(x(u,v),\,y(u,v)\bigr)
$$
be a local inverse of $F$. Then
$G\circ F=Id_{\mathbb{R}^2}$, that is,
$G\bigl(F(x,y)\bigr)=(x,y)$.

Differentiating the last identity and applying the chain rule, we obtain
$$
D(G\circ F)
=
DG\bigl(F(x,y)\bigr)\,DF(x,y).
$$
Since
$$
D(Id_{\mathbb{R}^2})=I,
$$
where $I$ is the identity matrix of order two, it follows that
$$
DG\bigl(F(x,y)\bigr)\,DF(x,y)=I.
$$
Hence,
$$
DG\bigl(F(x,y)\bigr)
=
\Biggr(DF(x,y)\Biggr)^{-1}.
$$
The lemma is proved.

\vskip2mm
\textbf{Lemma 2.}
Let
$$
u=u(x,y),\qquad
v=v(x,y),
$$
and suppose that
$$
J=\frac{\partial(u,v)}{\partial(x,y)}=u_xv_y-u_yv_x\neq0.
$$
Then
\begin{equation}\label{hodograf-nelin-x-y}
x_u=\frac{v_y}{J},
\qquad
x_v=-\frac{u_y}{J},
\qquad
y_u=-\frac{v_x}{J},
\qquad
y_v=\frac{u_x}{J}.
\end{equation}

\textbf{Proof.}
By Lemma 1,
%\ref{Lema-1-nelin}
$$
DG(u,v)=\Biggr(DF(x,y)\Biggr)^{-1},
$$
where
$$
DF(x,y)=\left(
\begin{array}{cc}
u_x&u_y\\
v_x&v_y
\end{array}
\right),
\qquad
DG(u,v)=
\left(
\begin{array}{cc}
x_u&x_v\\
y_u&y_v
\end{array}
\right).
$$
Therefore, using the formula for the inverse of a second-order matrix, we obtain
$$
\left(
\begin{array}{cc}
x_u&x_v\\
y_u&y_v
\end{array}
\right)
=
\frac{1}{J}
\left(
\begin{array}{cc}
v_y&-u_y\\
-v_x&u_x
\end{array}
\right),
$$
from which the required identities follow immediately.
The lemma is proved.

\vskip2mm
\textbf{Theorem 3.}\label{Teor-3-nelin}
Let the function
$\Phi(\zeta)=u(x,y)e_1+v(x,y)e_2$
be monogenic in the domain $D$, and suppose that its Jacobian
$J=\frac{\partial(u,v)}{\partial(x,y)}$
does not vanish in $D$.
Then, under the hodograph transformation (\ref{hodograf-nelin}), the Cauchy--Riemann system
(\ref{A-2-zmin-koef-systema-K-R-nelin})
is transformed into the linear system
\begin{equation}\label{A-2-zmin-koef-godograf-systema}
x_v=A(u)y_u,
\qquad
x_u=y_v-B(u)y_u.
\end{equation}

\textbf{Proof.}
Using the Cauchy--Riemann system
(\ref{A-2-zmin-koef-systema-K-R-nelin}), from
(\ref{hodograf-nelin-x-y}) we obtain
$$
x_u=\frac{u_x+B(u)v_x}{J},
\qquad
x_v=-\frac{A(u)v_x}{J},
\qquad
y_u=-\frac{v_x}{J},
\qquad
y_v=\frac{u_x}{J}.
$$
Hence, system (\ref{A-2-zmin-koef-godograf-systema}) follows.
The theorem is proved.

We now reduce system (\ref{A-2-zmin-koef-godograf-systema}) to a linear PDE with variable coefficients.

\vskip2mm
\textbf{Theorem 4.}\label{Teor-4-nelin}
Let the assumptions of Theorem 3 be satisfied.
%\ref{Teor-3-nelin}
Then the second component $y=y(u,v)$ of the hodograph transformation is a solution of the linear differential equation
\begin{equation}\label{rivn-dlia-y-pislia-godografa}
y_{vv}=A(u)y_{uu}+B(u)y_{uv}+A'(u)y_u.
\end{equation}

\textbf{Proof.}
Differentiating the first equation of system
(\ref{A-2-zmin-koef-godograf-systema})
with respect to $u$, we obtain
$$
x_{vu}=A'(u)y_u+A(u)y_{uu}.
$$
Since the function $A$ depends only on $u$, we have
$$
\frac{\partial A}{\partial v}=0.
$$

Next, differentiating the second equation of system
(\ref{A-2-zmin-koef-godograf-systema})
with respect to $v$, we obtain
$$
x_{uv}=y_{vv}-B(u)y_{uv},
$$
since
$$
\frac{\partial B}{\partial v}=0.
$$

Using the equality of mixed derivatives $x_{uv}=x_{vu}$, we obtain
$$
A'(u)y_u+A(u)y_{uu}=y_{vv}-B(u)y_{uv},
$$
or
$$
y_{vv}=A(u)y_{uu}+B(u)y_{uv}+A'(u)y_u.
$$
The theorem is proved.

Thus, the problem of constructing solutions of the nonlinear equation
(\ref{rivn-dlia-u-A-B-vid-u-neliniyne})
is reduced to the construction of solutions of the linear equation
(\ref{rivn-dlia-y-pislia-godografa}).

\vskip2mm
\textbf{Theorem 5.}\label{Teor-5-nelin}
Let $\Omega\subset\mathbb{R}^2_{(u,v)}$ be a simply connected domain,
$A,B\in C^1(I)$, and let the function
$y=y(u,v)\in C^2(\Omega)$ satisfy the linear differential equation
\begin{equation}\label{rivn-dlia-y-pislia-godografa-obern}
y_{vv}
=
A(u)y_{uu}
+
B(u)y_{uv}
+
A'(u)y_u.
\end{equation}
Then there exists a function $x=x(u,v)\in C^2(\Omega)$ satisfying system
(\ref{A-2-zmin-koef-godograf-systema}).
If the Jacobian of the mapping
$(u,v)\longmapsto\bigl(x(u,v),y(u,v)\bigr)$
does not vanish, that is,
\begin{equation}\label{yakobian-obernyi-godograf}
K=\frac{\partial(x,y)}{\partial(u,v)}
=
x_uy_v-x_vy_u\neq0,
\end{equation}
then there exists an inverse hodograph transformation
\begin{equation}\label{obernyi-godograf}
u=u(x,y),\qquad v=v(x,y),
\end{equation}
and the functions $u$ and $v$ satisfy the nonlinear Cauchy--Riemann system
(\ref{A-2-zmin-koef-systema-K-R-nelin}).

\textbf{Proof.}
Consider system (\ref{A-2-zmin-koef-godograf-systema}), that is,
$$
x_v=A(u)y_u,
\qquad
x_u=y_v-B(u)y_u.
$$
For the existence of a function $x=x(u,v)$, it is necessary to verify the compatibility condition
$x_{uv}=x_{vu}$.

Differentiating the first equation of system
(\ref{A-2-zmin-koef-godograf-systema})
with respect to $u$, we obtain
$$
x_{vu}=A'(u)y_u+A(u)y_{uu}.
$$
Differentiating the second equation of this system with respect to $v$, we obtain
$$
x_{uv}=y_{vv}-B(u)y_{uv},
$$
since $B$ depends only on $u$.

Equation (\ref{rivn-dlia-y-pislia-godografa-obern}) implies that
$$
y_{vv}-B(u)y_{uv}
=
A(u)y_{uu}+A'(u)y_u.
$$
Hence,
$$
x_{uv}=x_{vu}.
$$

Thus, system (\ref{A-2-zmin-koef-godograf-systema}) is compatible. Therefore, there exists a function
$x=x(u,v)\in C^2(\Omega)$ satisfying this system. Since $\Omega$ is simply connected, the function $x$ can be defined throughout $\Omega$ up to an additive constant.

Suppose now that condition (\ref{yakobian-obernyi-godograf}) holds. Then, by the inverse function theorem, the mapping
$$
(u,v)\longmapsto\bigl(x(u,v),y(u,v)\bigr)
$$
is locally invertible. Hence, there exist functions
$$
u=u(x,y),\qquad v=v(x,y).
$$

Using the formulas for the derivatives of the inverse mapping, we have
$$
u_x=\frac{y_v}{K},
\qquad
u_y=-\frac{x_v}{K},
\qquad
v_x=-\frac{y_u}{K},
\qquad
v_y=\frac{x_u}{K}.
$$
Using the first equation of system
(\ref{A-2-zmin-koef-godograf-systema}), we obtain
$$
u_y=A(u)v_x.
$$
Next, using the second equation of system
(\ref{A-2-zmin-koef-godograf-systema}), we obtain
$$
v_y=u_x+B(u)v_x.
$$
Thus, the functions $u=u(x,y)$ and $v=v(x,y)$ satisfy the Cauchy--Riemann system
(\ref{A-2-zmin-koef-systema-K-R-nelin}).
The theorem is proved.

\vskip2mm
\textbf{Corollary 1.}
Theorems 4 and 5
%\ref{Teor-4-nelin} and \ref{Teor-5-nelin}
establish a local one-to-one correspondence between solutions of the nonlinear Cauchy--Riemann system
(\ref{A-2-zmin-koef-systema-K-R-nelin})
and solutions of the linear differential equation
(\ref{rivn-dlia-y-pislia-godografa-obern}).
Thus, the construction of solutions of the nonlinear system is reduced to solving the linear equation followed by the hodograph and inverse hodograph transformations.

Solutions of the nonlinear equation
(\ref{rivn-dlia-u-A-B-vid-u-neliniyne})
can therefore be constructed according to the following scheme.

\vskip2mm
\textbf{Algorithm for constructing solutions.}

\begin{enumerate}
\item Choose the functions $A(u)$ and $B(u)$.

\item Solve the linear differential equation
(\ref{rivn-dlia-y-pislia-godografa}).

\item Recover the function $x=x(u,v)$ from system
(\ref{A-2-zmin-koef-godograf-systema}).

\item Perform the inverse hodograph transformation
(\ref{obernyi-godograf}).

\item The resulting function $u=u(x,y)$ is a solution of the nonlinear differential equation
(\ref{rivn-dlia-u-A-B-vid-u-neliniyne}).
\end{enumerate}

\section{Example 1}

Let
$$
A(u)=u,\qquad B(u)=0.
$$
Then the multiplication table of the algebra is determined by
$e_2^2=ue_1$,
and the Cauchy--Riemann system takes the form
\begin{equation}\label{systema-K-R-pryklad-1-nelin}
u_y=uv_x,\qquad v_y=u_x.
\end{equation}

The corresponding nonlinear differential equation
(\ref{rivn-dlia-u-A-B-vid-u-neliniyne}) for the first component has the form
\begin{equation}\label{neliniyne-rivn-pryklad-1}
uu_{yy}=u^2u_{xx}+u_y^2.
\end{equation}
Under the condition $u\neq0$, equation
(\ref{neliniyne-rivn-pryklad-1}) can be written in the normalized form
$$
u_{yy}=uu_{xx}+\frac{u_y^2}{u}.
$$

After the hodograph transformation, the linear equation
(\ref{rivn-dlia-y-pislia-godografa-obern}) takes the form
\begin{equation}\label{liniyne-rivn-pryklad-1}
y_{vv}=uy_{uu}+y_u.
\end{equation}

\vskip2mm
\textbf{Case 1.1.}
Consider the following solution of the linear equation
(\ref{liniyne-rivn-pryklad-1}):
$$
y(u,v)=\ln u,\qquad u>0.
$$
Since
$$
y_{vv}=0,\qquad
y_u=\frac{1}{u},\qquad
y_{uu}=-\frac{1}{u^2},
$$
we have
$$
uy_{uu}+y_u=-\frac{1}{u}+\frac{1}{u}=0.
$$
Hence, the function $y(u,v)=\ln u$ is a solution of equation
(\ref{liniyne-rivn-pryklad-1}).

Let us recover the function $x=x(u,v)$ from system
(\ref{A-2-zmin-koef-godograf-systema}). We have
$$
x_v=uy_u=1,\qquad
x_u=y_v=0.
$$
Therefore,
$x(u,v)=v+C$.
Setting $C=0$, we obtain
$$
x=v,\qquad y=\ln u.
$$

The Jacobian of this mapping is
$$
K=\frac{\partial(x,y)}{\partial(u,v)}
=x_uy_v-x_vy_u
=-\frac{1}{u}\neq0,
\qquad u>0.
$$
Hence, the inverse hodograph transformation exists and is given by
$$
u=e^y,\qquad v=x.
$$
Thus, the function
$u(x,y)=e^y$
is a solution of the nonlinear equation
(\ref{neliniyne-rivn-pryklad-1}), while the functions
$$
u(x,y)=e^y,\qquad v(x,y)=x
$$
satisfy the Cauchy--Riemann system
(\ref{systema-K-R-pryklad-1-nelin}).

\vskip2mm
\textbf{Case 1.2.}
Consider another solution of the linear equation:
$$
y(u,v)=v^2+2u.
$$

Let us recover the function $x=x(u,v)$ from the system
$$
x_v=uy_u,\qquad x_u=y_v.
$$
We have
$$
x_v=2u,\qquad x_u=2v,
$$
and hence
$x(u,v)=2uv+C$.
Setting $C=0$, we obtain
$$
x=2uv,\qquad y=v^2+2u.
$$

The Jacobian of this mapping is
$$
K=x_uy_v-x_vy_u=4(v^2-u).
$$
Therefore, the hodograph transformation is locally invertible in those regions where
$v^2\neq u$.

From the equality $x=2uv$, we obtain
$u=\frac{x}{2v}$.
Substituting this expression into $y=v^2+2u$, we obtain
$$
v^3-yv+x=0.
$$
Hence, every root $v=v(x,y)$ of this cubic equation for which
$v\neq0$ and $v^2\neq u$
determines a solution
$$
u(x,y)=\frac{x}{2v(x,y)}
$$
of the nonlinear equation
(\ref{neliniyne-rivn-pryklad-1}).

Thus, every solution of the linear equation
(\ref{liniyne-rivn-pryklad-1})
generates a solution of the nonlinear equation
(\ref{neliniyne-rivn-pryklad-1}).

\section{Example 2}

Let
$$
A(u)=1,\qquad B(u)=u.
$$
Then $A'(u)=0$, $B'(u)=1$, and the nonlinear equation
(\ref{rivn-dlia-u-A-B-vid-u-neliniyne-normal})
takes the form
\begin{equation}\label{neliniyne-rivn-pryklad-2}
u_{yy}=u_{xx}+u\,u_{xy}+u_xu_y.
\end{equation}

After the hodograph transformation, the corresponding linear equation for the function
$y=y(u,v)$ has the form
\begin{equation}\label{liniyne-rivn-pryklad-2}
y_{vv}=y_{uu}+u\,y_{uv}.
\end{equation}

\vskip2mm
\textbf{Case 2.1.}
Consider the function
$$
y(u,v)=u^2+v^2.
$$
For this function,
$$
y_{vv}=2,\qquad
y_{uu}=2,\qquad
y_{uv}=0.
$$
Hence,
$$
y_{uu}+u\,y_{uv}=2=y_{vv},
$$
so the function $y(u,v)=u^2+v^2$ is a solution of equation
(\ref{liniyne-rivn-pryklad-2}).

Let us recover the function $x=x(u,v)$ from system
(\ref{A-2-zmin-koef-godograf-systema}). We have
$$
x_v=A(u)y_u=2u
$$
and
$$
x_u=y_v-B(u)y_u=2v-2u^2.
$$
Integrating the first equality with respect to $v$, we obtain
$$
x(u,v)=2uv+\varphi(u).
$$
Therefore,
$$
x_u=2v+\varphi'(u)=2v-2u^2,
$$
which gives
$\varphi'(u)=-2u^2$.

Hence,
$$
\varphi(u)=-\frac{2}{3}u^3+C.
$$
Setting $C=0$, we obtain the hodograph mapping
\begin{equation}\label{godograf-pryklad-2}
x=2uv-\frac{2}{3}u^3,
\qquad
y=u^2+v^2.
\end{equation}

The Jacobian of this mapping is
$$
K=\frac{\partial(x,y)}{\partial(u,v)}
=
x_uy_v-x_vy_u
=
(2v-2u^2)2v-(2u)2u
=
4\bigl(v^2-u^2v-u^2\bigr).
$$
Hence, the mapping
(\ref{godograf-pryklad-2})
is locally invertible in those regions where
$v^2-u^2v-u^2\neq0$.

After inversion of the hodograph transformation, the function
$u=u(x,y)$ is a solution of the nonlinear equation
(\ref{neliniyne-rivn-pryklad-2}).
From
(\ref{godograf-pryklad-2}), the function $u=u(x,y)$ is determined implicitly by
$$
4u^{2}y
=
4u^{4}
+
\left(x+\frac{2}{3}u^{3}\right)^{2}.
$$

\vskip2mm
\textbf{Case 2.2.}
Consider the function
$$
y(u,v)=u^3-6uv.
$$
We have
$$
y_{vv}=0,
\qquad
y_{uu}=6u,
\qquad
y_{uv}=-6.
$$
Therefore,
$$
y_{uu}+u\,y_{uv}
=
6u-6u=0=y_{vv},
$$
so the function
$y(u,v)=u^3-6uv$
is a solution of equation
(\ref{liniyne-rivn-pryklad-2}).

Let us recover the function $x=x(u,v)$ from system
(\ref{A-2-zmin-koef-godograf-systema}).
Since $A(u)=1$ and $B(u)=u$, we have
$$
x_v=A(u)y_u=3u^2-6v,
$$
and
$$
x_u
=
y_v-B(u)y_u
=
6uv-3u^3-6u.
$$

Integrating the first equality with respect to $v$, we obtain
$$
x(u,v)=3u^2v-3v^2+\varphi(u).
$$
Therefore,
$$
x_u=6uv+\varphi'(u)
=
6uv-3u^3-6u,
$$
which gives
$$
\varphi'(u)=-3u^3-6u.
$$
Hence,
$$
\varphi(u)=-\frac34u^4-3u^2+C.
$$

Setting $C=0$, we obtain the hodograph mapping
\begin{equation}\label{godograf-pryklad-4}
x=3u^2v-3v^2-\frac34u^4-3u^2,
\qquad
y=u^3-6uv.
\end{equation}

The Jacobian of this mapping is
$$
K=\frac{\partial(x,y)}{\partial(u,v)}
=
x_uy_v-x_vy_u.
$$
Since
$$
x_u=6uv-3u^3-6u,
\qquad
x_v=3u^2-6v,
$$
and
$$
y_u=3u^2-6v,
\qquad
y_v=-6u,
$$
we obtain
$$
K
=
(6uv-3u^3-6u)(-6u)
-
(3u^2-6v)^2
=
9\bigl(u^4+4u^2-4v^2\bigr).
$$

Hence, the hodograph mapping
(\ref{godograf-pryklad-4})
is locally invertible in those regions where
$$
u^4+4u^2-4v^2\neq0.
$$

From the second equality in
(\ref{godograf-pryklad-4}), we obtain
$$
v=\frac{u^3-y}{6u},
\qquad
u\neq0.
$$
Substituting this expression into the formula for $x$, we obtain
$$
x
=
-\frac{4u^6+36u^4+4u^3y+y^2}{12u^2}.
$$

Therefore, after inversion of the hodograph transformation, the function
$u=u(x,y)$ is a local solution of the nonlinear equation
(\ref{neliniyne-rivn-pryklad-2}),
implicitly determined by the equation
$$
12u^2x+4u^6+36u^4+4u^3y+y^2=0.
$$

\section{Example 3}

Let
$$
A(u)=1,\qquad B(u)=u^2.
$$
Then $A'(u)=0$, $B'(u)=2u$, and the nonlinear equation
(\ref{rivn-dlia-u-A-B-vid-u-neliniyne-normal})
takes the form
\begin{equation}\label{neliniyne-rivn-pryklad-3}
u_{yy}
=
u_{xx}
+
u^2u_{xy}
+
2u\,u_xu_y.
\end{equation}

After the hodograph transformation, the corresponding linear equation for the function $y=y(u,v)$ has the form
\begin{equation}\label{liniyne-rivn-pryklad-3}
y_{vv}
=
y_{uu}
+
u^2y_{uv}.
\end{equation}

It is easy to verify that the function
$$
y(u,v)=uv-\frac{1}{12}u^4
$$
is a solution of equation (\ref{liniyne-rivn-pryklad-3}).

Let us recover the function $x=x(u,v)$ from system
(\ref{A-2-zmin-koef-godograf-systema}).
We have
$$
x_v=A(u)y_u=v-\frac{1}{3}u^3
$$
and
$$
x_u
=
y_v-B(u)y_u
=
u-u^2\left(v-\frac{1}{3}u^3\right)
=
u-u^2v+\frac{1}{3}u^5.
$$

Integrating the first equality with respect to $v$, we obtain
$$
x(u,v)
=
\frac{1}{2}v^2
-
\frac{1}{3}u^3v
+
\varphi(u).
$$
Hence,
$$
x_u=-u^2v+\varphi'(u).
$$
Comparing this expression with
$$
x_u=u-u^2v+\frac{1}{3}u^5,
$$
we obtain
$$
\varphi'(u)=u+\frac{1}{3}u^5.
$$
Therefore,
$$
\varphi(u)=\frac{1}{2}u^2+\frac{1}{18}u^6+C.
$$

Setting $C=0$, we obtain the hodograph mapping
\begin{equation}\label{godograf-pryklad-3}
x=
\frac{1}{2}v^2
-
\frac{1}{3}u^3v
+
\frac{1}{2}u^2
+
\frac{1}{18}u^6,
\qquad
y=uv-\frac{1}{12}u^4.
\end{equation}

The Jacobian of this mapping is
$$
K=x_uy_v-x_vy_u
=
\left(u-u^2v+\frac{1}{3}u^5\right)u
-
\left(v-\frac{1}{3}u^3\right)^2
=
u^2-v^2-\frac{1}{3}u^3v+\frac{2}{9}u^6.
$$

Hence, the hodograph mapping
(\ref{godograf-pryklad-3})
is locally invertible in those regions where
$$
u^2-v^2-\frac{1}{3}u^3v+\frac{2}{9}u^6
\neq0.
$$

From the second equality in
(\ref{godograf-pryklad-3}), under the condition $u\neq0$, we obtain
$$
v=\frac{y}{u}+\frac{1}{12}u^3.
$$
Substituting this expression into the first equality in
(\ref{godograf-pryklad-3}), we obtain
$$
x=
\frac{u^8-8u^4y+16u^4+16y^2}{32u^2}.
$$

Thus, after inversion of the hodograph transformation, the function
$u=u(x,y)$ is a local solution of the nonlinear equation
(\ref{neliniyne-rivn-pryklad-3}), which, under the condition $u\neq0$, is implicitly determined by the equation
$$
u^8-8u^4y+16u^4+16y^2-32u^2x=0.
$$

\section{Concluding Remarks}

1. A natural direction for further research is to study generalized hypercomplex systems in which the structural coefficients $A$ and $B$
depend not only on the function $u$ but also on its derivatives $u_x$ and $u_y$.
Such a generalization is expected to lead to new classes of fully nonlinear differential equations.

2. We note that, for simplicity, in this paper we have considered the case in which the structural functions $A$ and $B$ depend only on the component $u$.
However, in a similar way one can investigate much broader classes of hypercomplex systems in which $A$ and $B$ depend on the component $v$,
on the combination $u+v$, on the expressions $u^2+v^2$, $uv$, or, more generally, on an arbitrary function $F(u,v)$.
Each such choice leads to its own nonlinear Cauchy--Riemann system and, after elimination of one of the components, generates a new class of nonlinear differential equations.

\subsection*{Acknowledgment}
 This work was supported by a grant from the Simons Foundation 
(SFI-PD-Ukraine-00014586,V.S.Sh.). The author was supported by budget program “Support for the development of priority areas of research” 
(KPKVK 6541230).

\renewcommand{\refname}{References}


\begin{thebibliography}{99}
%\bibitem{Shpakivskyi-EJDE}
%Shpakivskyi V. S. \emph{Construction of solutions of PDEs using holomorphic functions of several variables}, 
%Electronic Journal of Differential Equations, Vol. 2024 (2024), No. 71, pp. 1--12. 

\bibitem{Pierce}
Pierce R. S. \emph{Associative Algebras}, Graduate Texts in Mathematics, Vol. 88, Springer-Verlag, New York, 1982. DOI: 10.1007/978-1-4757-0163-0.

\bibitem{Plaksa-Shpakivskyi-2023}
Plaksa S. A., Shpakivskyi V. S. \emph{Monogenic Functions in Spaces with Commutative Multiplication and Applications}.  Frontiers in Mathematics, Birk\-h\"{a}user Cham, 2023, 550 p.

\bibitem{Ames}
Ames W. F. \emph{Nonlinear Partial Differential Equations in Engineering}. Academic Press, New York, 1965.

\bibitem{Clarkson}
Clarkson P. A., Fokas A. S., Ablowitz M. J. \emph{Hodograph transformations of linearizable partial differential equations}. SIAM Journal on Applied Mathematics, 49(4) (1989), 1188–1209. DOI: 10.1137/0149071.

\bibitem{Qu}
Qu C., Kang J. \emph{Hodograph transformations and Cauchy problem to systems of nonlinear parabolic equations}. Studies in Applied Mathematics, 132(2) (2014), 81–111. DOI: 10.1111/sapm.12025.


\end{thebibliography}
\end{document}